\documentclass[10pt]{amsart}

\usepackage[latin1]{inputenc}
\usepackage{amsmath}
\usepackage{amsfonts}
\usepackage{amssymb}
\usepackage[mathscr]{eucal}
\usepackage{diagrams}
\usepackage{xy,color,graphicx}
\usepackage{hyperref}
\xyoption{all}

\newcommand{\wis}[1]{{\text{\em \usefont{OT1}{cmtt}{m}{n} #1}}}

\newcommand{\N}{\mathbb{N}}
\newcommand{\C}{\mathbb{C}}
\newcommand{\Z}{\mathbb{Z}}
\newcommand{\vtx}[1]{*+[o][F-]{\scriptscriptstyle #1}}
\newcommand{\Oscr}{\mathcal{O}}

\newtheorem{proposition}{Proposition}

\newtheorem{theorem}{Theorem}
\newtheorem{lemma}{Lemma}
\newtheorem{example}{Example}
\newtheorem{remark}{Remark}

\title{Matrix transposition and braid reversion}
\author{Lieven Le Bruyn} 
\address{Department of Mathematics, University of Antwerp \\ 
 Middelheimlaan 1, B-2020 Antwerp (Belgium) \\ {\tt lieven.lebruyn@ua.ac.be}}

\begin{document}
\sloppy

\maketitle

\begin{abstract} Matrix transposition induces an involution $\tau$ on the  isomorphism classes of semi-simple $n$-dimensional representations of the three string braid group $B_3$. We show that a  connected component of this variety can detect braid-reversion or that $\tau$ acts as the identity on it. We classify the fixed-point components. \end{abstract}

\section{Introduction}

If $\phi=(X_1,X_2)$ is an $n$-dimensional representation of the three string braid group $B_3 = \langle \sigma_1,\sigma_2~|~\sigma_1 \sigma_2 \sigma_1 = \sigma_2 \sigma_1 \sigma_2 \rangle$, then so is the pair of transposed matrices $\tau(\phi) = (X_1^{tr},X_2^{tr})$. In this paper we investigate when $\phi \simeq \tau(\phi)$.

This problem is relevant to detect braid- and knot-reversion. Recall that  a knot is said to be invertible if it can be deformed continuously to itself, but with the orientation reversed. There do exist non-invertible knots, the unique one with a minimal number of crossings is knot $8_{17}$, see the Knot Atlas \cite{KnotAtlas},  which is the closure of the three string braid $b=\sigma_1^{-2} \sigma_2 \sigma_1^{-1} \sigma_2 \sigma_1^{-1} \sigma_2^2$.  Proving knot-vertibility of $8_{17}$ comes down to separating the conjugacy class of the braid $b$ from that of its reversed braid $b' = \sigma_2^2 \sigma_1^{-1} \sigma_2 \sigma_1^{-1} \sigma_2 \sigma_1^{-2}$. Now, observe that a $B_3$-representation $\phi$ can separate $b$ from $b'$, via $Tr_{\phi}(b) \not= Tr_{\phi}(b')$, only if $\phi$ is not isomorphic to the transposed representation $\tau(\phi)$.

The involution $\tau$ on the affine variety $\wis{rep}_n B_3$ of all $n$-dimensional representations passes to an involution $\tau$ on the quotient variety $\wis{rep}_n B_3/PGL_n = \wis{iss}_n B_3$, classifying isomorphism classes of semi-simple $n$-dimensional representations. Recall from \cite{LeBruynB3} and \cite{Westbury} that $\wis{iss}_n B_3$ decomposes as a disjoint union of its irreducible components
\[
\wis{iss}_n B_3 = \sqcup_{\alpha}~\wis{iss}_{\alpha} B_3 \]
and the components containing a Zariski open subset of simple $B_3$-representations are classified by $\alpha = (a,b;x,y,z) \in \N^5$ satisfying $n=a+b=x+y+z$, $a \geq b$ and $x=max(x,y,z) \leq b$. 

\begin{theorem} \label{main} If the component $\wis{iss}_{\alpha} B_3$ contains simple representations, then $\tau$ acts as the identity on it  if and only if $\alpha$ is equal to
\begin{itemize}
\item{$(1,0;1,0,0)$, or $(4,2;2,2,2)$, or}
\item{$(k,k;k,k-1,1)$, or $(k,k;k,1,k-1)$, or}
\item{$(k+1,k;k,k,1)$, or $(k+1,k;k,1,k)$}
\end{itemize}
for  some $k \geq 1$. In all these cases, $dim~\wis{iss}_{\alpha} B_3 = n$. In all other cases, $\wis{iss}_{\alpha} B_3$ contains simple representations $\phi$ such that $Tr_{\phi}(b) \not= Tr_{\phi}(b')$, that is, $\wis{iss}_{\alpha} B_3$ can detect braid-reversion.
\end{theorem}

Note that this result generalizes \cite{LeBruynB3} where it was proved that there is a unique component of $\wis{iss}_6 B_3$, namely $\wis{iss}_{\alpha} B_3$ for $\alpha=(3,3;2,2,2)$, containing representations $\phi$ such that $Tr_{\phi}(b) \not= Tr_{\phi}(b')$.

\section{The involution $\tau$ and stable quiver representations}

In this section we follow Bruce Westbury \cite{Westbury} reducing the study of simple $B_3$-representations to specific stable quiver-representations, and, we will describe the involution $\tau$ in terms of these representations.

If $\phi = (X_1,X_2)$ is a simple $n$-dimensional $B_3$-representation, then the central element $c=(\sigma_1 \sigma_2)^3 = (\sigma_1 \sigma_2 \sigma_1)^2$ acts via a scalar matrix $\lambda I_n$ for some $\lambda \not= 0$. Hence, $\phi'=\lambda^{-1/6} \phi = (X'_1,X'_2)$ is a simple representation of the quotient group
\[
B_3/ \langle c \rangle = \langle s,t~|~s^2=t^3=e \rangle \simeq C_2 \ast C_3 \simeq \Gamma  \]
which is the free product of cyclic groups of order two and three (and thus isomorphic the modular group $\Gamma = PSL_2(\Z)$) where $s$ and $t$ are the images of $\sigma_1 \sigma_2 \sigma_1$ and $\sigma_1 \sigma_2$. Decompose the $n$-dimensional space $V = \C^n_{\phi'}$ into eigenspaces for the actions of $s$ and $t$
\[
V_{+} \oplus V_{-} = V = V_1 \oplus V_{\rho} \oplus V_{\rho^2} \]
with $\rho$ a primitive 3rd root of unity. For $a=dim(V_+), b=dim(V_-)$, $x=dim(V_1), y=dim(V_{\rho})$ and $ z = dim(V_{\rho^2})$, clearly $a+b=n=x+y+z$. 
Choose a vector-space basis for $V$ compatible with the decomposition $V_+ \oplus V_-$ and another basis  of $V$ compatible with the decomposition $V_1 \oplus V_{\rho} \oplus V_{\rho^2}$, then the corresponding base-change matrix 
\[
B = \begin{bmatrix} B_{11} & B_{12} \\ B_{21} & B_{22} \\ B_{31} & B_{32} \end{bmatrix}
\in GL_n(\C) \]
 determines the quiver-representation $V_B$ with dimension vector $\alpha = (a,b;x,y,z)$
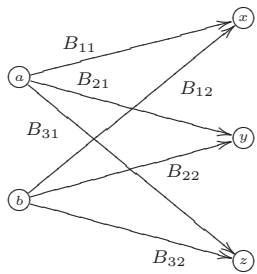
\begin{eqnarray*} \label{identification}
\xymatrix@=.5cm{
& & & & \vtx{x} \\
\vtx{a} \ar[rrrru]^(.3){B_{11}} \ar[rrrrd]^(.3){B_{21}} \ar[rrrrddd]_(.2){B_{31}} & & & & \\
& & & & \vtx{y} \\
\vtx{b} \ar[rrrruuu]_(.7){B_{12}} \ar[rrrru]_(.7){B_{22}} \ar[rrrrd]_(.7){B_{32}} & & & & \\
& & & & \vtx{z}}
\end{eqnarray*}
For $B$ invertible, the representation $V_B$ is semi-stable in the sense of \cite{King}, meaning that for every proper sub-representations $W$, with dimension vector $\beta = (a',b';x',y',z')$ we have $x'+y'+z' \geq a'+b'$. If all these inequalities are strict, we call $V_B$ a stable representation, which is equivalent to the $\Gamma$-representation $V=\C^n_{\phi'}$ being simple. Westbury \cite{Westbury} showed that two $\Gamma$-representations $V \simeq V'$ if and only if $V_B \simeq V'_{B'}$ as quiver-representations, that is, there exist base-changes in the eigenspaces
\[
(M_1,M_2;N_1,N_2,N_3) \in GL_a \times GL_b \times GL_x \times GL_y \times GL_z \]
such that
\[
\begin{bmatrix} N_1 & 0 & 0 \\ 0 & N_2 & 0 \\ 0 & 0 & N_3 \end{bmatrix}  \begin{bmatrix} B'_{11} & B'_{12} \\ B'_{21} & B'_{22} \\ B'_{31} & B'_{32} \end{bmatrix} \begin{bmatrix} M_1^{-1} & 0 \\ 0 & M_2^{-1} \end{bmatrix}  = \begin{bmatrix} B_{11} & B_{12} \\ B_{21} & B_{22} \\ B_{31} & B_{32} \end{bmatrix} \]
 Working backwards, we recover the $B_3$-representation $\phi = (X_1,X_2)$ from the invertible matrix $B$ via
\[
( \ast )  \begin{cases}
X_1 = \lambda^{1/6} B^{-1} \begin{bmatrix} 1_{x} & 0 & 0 \\ 0 & \rho^2 1_{y} & 0 \\ 0 & 0 & \rho 1_{z} \end{bmatrix} B \begin{bmatrix} 1_{a} & 0 \\ 0 & -1_{b} \end{bmatrix} \\
X_2 = \lambda^{1/6} \begin{bmatrix} 1_{a} & 0 \\ 0 & -1_{b} \end{bmatrix} B^{-1}  \begin{bmatrix} 1_{x} & 0 & 0 \\ 0 & \rho^2 1_{y} & 0 \\ 0 & 0 & \rho 1_{z} \end{bmatrix} B
\end{cases}
\]

\begin{proposition} If the $n$-dimensional simple $B_3$-representation $\phi = (X_1,X_2)$ is determined by $\lambda \in \C^*$ and the stable quiver-representation $V_B$, then $\tau(\phi)=(X_1^{tr},X_2^{tr})$ is isomorphic to the representation determined by $\lambda$ and the stable quiver-representation $V_{(B^{-1})^{tr}}$.
\end{proposition}

\begin{proof} Taking transposes of the formulas $(\ast)$ for the $X_i$ we get
\[
\begin{cases}
X_1^{tr} = \lambda^{1/6}  \begin{bmatrix} 1_{a} & 0 \\ 0 & -1_{b} \end{bmatrix} B^{tr} \begin{bmatrix} 1_{x} & 0 & 0 \\ 0 & \rho^2 1_{y} & 0 \\ 0 & 0 & \rho 1_{z} \end{bmatrix} (B^{-1})^{tr}  \\
X_2^{tr} = \lambda^{1/6}  B^{tr}  \begin{bmatrix} 1_{x} & 0 & 0 \\ 0 & \rho^2 1_{y} & 0 \\ 0 & 0 & \rho 1_{z} \end{bmatrix} (B^{-1})^{tr} \begin{bmatrix} 1_{a} & 0 \\ 0 & -1_{b} \end{bmatrix}
\end{cases}
\]
Conjugating these with the matrix $\begin{bmatrix} 1_{a} & 0 \\ 0 & -1_{b} \end{bmatrix}$ (which is also a base-change action in $GL(\alpha)$) we obtain again a matrix-pair in standard-form $(\ast)$, this time replacing the matrix $B$ by the matrix $(B^{-1})^{tr}$.
\end{proof}

That is, we have reduced the original problem of verifying whether or not $\phi \simeq \tau(\phi)$ as $B_3$-representations to the problem of verifying whether or not the two stable representations $V_B$ and $V_{(B^{-1})^{tr}}$ lie in the same $GL(\alpha)$-orbit.

\begin{example} \label{example} The two components $\wis{iss}_{\alpha} B_3$ containing simple $2$-dimensional $B_3$-representations for $\alpha=(1,1;1,1,0)$ or $(1,1;1,0,1)$ are fixed-point components for the involution $\tau$.
A general stable $\alpha=(1,1;1,1,0)$  dimensional representation $V_B$ is isomorphic to one of the form
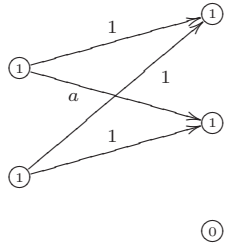
\[
 \xymatrix@=.4cm{
& & & & \vtx{1} \\
\vtx{1} \ar[rrrru]^1 \ar[rrrrd]_(0.3){a}  & & & & \\
& & & & \vtx{1} \\
\vtx{1} \ar[rrrru]^{1} \ar[rrrruuu]_(0.7){1}  & & & & \\
& & & &  \vtx{0} }
\]
with $a \not= 1$. Hence, we can take $B=\begin{bmatrix} 1 & 1 \\ a & 1 \end{bmatrix}$. But then, $V_B$ and $V_{(B^{-1})^{tr}}$ lie in the same $GL(\alpha) = \C^* \times \C^* \times \C^* \times \C^*$-orbit because
\[
(B^{-1})^{tr} =  \frac{1}{1-a} \begin{bmatrix} 1 & -a \\ -1 & 1 \end{bmatrix} =  \begin{bmatrix} 1 & 0 \\ 0 & -\frac{1}{a} \end{bmatrix} B \begin{bmatrix} \frac{1}{1-a} & 0 \\ 0 & \frac{-a}{1-a} \end{bmatrix} \]
\end{example}

\section{The stratification and potential fixed-point components}

In this section we will show that a component $\wis{iss}_{\alpha} B_3$ containing $n$-dimensional simple $B_3$-representations is a fixed-point component for the involution $\tau$ only if $\alpha$ is among the list of theorem~\ref{main}. 

Because the group algebra $\C \Gamma = \C C_2 \ast \C C_3$ is a formally smooth algebra, we have a Luna stratification of $\wis{iss}_{\alpha} \Gamma$ by representation types, see \cite[\S 5.1]{LeBruynBook}. A point $p$ in $\wis{iss}_{\alpha} \Gamma$ determines the isomorphism class of a semi-simple $\Gamma$-representation
\[
V_p = S_1^{\oplus e_1} \oplus \hdots \oplus S_k^{\oplus e_k} \]
with all $S_i$ distinct simple $\Gamma$-representations with corresponding dimension vectors $\beta_i = (a_i,b_i;x_i,y_i,z_i)$. We say that $p$ (or $V_p$) is of representation type
\[
\tau = (e_1,\beta_1;\hdots;e_k,\beta_k) \quad \text{and clearly} \quad \alpha = \sum_i e_i \beta_i \]
With $\wis{iss}_{\alpha} \Gamma(\tau)$ we denote the subset of all points of representation type $\tau$. Recall that $\beta_i$ is the dimension vector of a simple $\Gamma$-representation if and only if $a_i+b_i=x_i+y_i+z_i$ and $max(x_i,y_i,z_i) \leq min(a_i,b_i)$ if $x_iy_iz_i \not= 0$ (the remaining cases being the $1$- and $2$-dimensional components). It follows from Luna's results \cite{Luna}  that every $\wis{iss}_{\alpha} \Gamma(\tau)$ is a locally closed smooth irreducible subvariety of $\wis{iss}_{\alpha} \Gamma$ of dimension $\sum_i (1+2a_ib_i-(x_i^2+y_i^2+z_i^2)$ and that
\[
\wis{iss}_{\alpha} \Gamma = \bigsqcup_{\tau} \wis{iss}_{\alpha} \Gamma(\tau) \]
is a finite smooth stratification of $\wis{iss}_{\alpha} \Gamma$. Degeneration of representation types, see \cite[p. 247]{LeBruynBook}, defines an ordering $\leq$ on representation types and by \cite[Prop. 5.3]{LeBruynBook} we have that $\wis{iss}_{\alpha} \Gamma(\tau')$ lies in the Zariski closure of $\wis{iss}_{\alpha} \Gamma(\tau)$ if and only if $\tau' \leq \tau$.

Observe that the involution $\tau$ on $\wis{iss}_{\alpha} \Gamma$ induced by $\tau(V_B) = V_{(B^{-1})^{tr}}$ preserves the strata and its restriction to $\wis{iss}_{\alpha} \Gamma(\tau)$ is induced by the involutions $\tau$ on the components $\wis{iss}_{\beta_i} \Gamma$. As the fixed-point set of $\tau$ is a closed subvariety of $\wis{iss}_{\alpha} \Gamma$ we deduce immediately :

\begin{lemma} \label{strata} If $\tau$ is the identity on a Zariski open subset of $\wis{iss}_{\alpha} \Gamma(\tau)$, then $\tau = id$ on all strata $\wis{iss}_{\alpha} \Gamma(\tau')$ with $\tau' \leq \tau$. Conversely, if $\tau = (e_1,\beta_1;\hdots;e_k,\beta_k)$ and $\tau \not= id$ on one of the components $\wis{iss}_{\beta_i} \Gamma$, then $\tau \not= id$ on all strata $\wis{iss}_{\alpha} \Gamma(\tau')$ with $\tau \leq \tau'$.
\end{lemma}

 In \cite{LeBruynB3} we have shown that for $\beta=(3,3;2,2,2)$ there are simple $B_3$-representations able to separate the braid $b$ from the introduction from its reversed braid $b'$. In particular, $\tau$ does not act as the identity on $\wis{iss}_{\beta} \Gamma$. We proved this by  parametrizing the matrices $B$ for a dense open subset of $\wis{iss}_{\beta} \Gamma$ by
 \[
 B= \begin{bmatrix} 1 & 0 & 0 & a & 0 & f \\
0 & 1 & 1 & 0 & 1 & 0 \\
1 & 1 & 0 & 1 & 0 & 0 \\
0 & 0 & 1 & 0 & d & e \\
0 & 1 & 0 & b & c & 0 \\
g & 0 & 1 & 0 & 0 & 1 \end{bmatrix}
\]
for free parameters $a,\hdots,g$. We then computed the matrix-pair $\phi=(X_1,X_2)$ from $(\ast)$ with generic values of the parameters in $\Z[\rho]$ and checked that $Tr_{\phi}(b) \not= Tr_{\phi}(b')$.

\begin{proposition} If $\alpha$ is the dimension vector of a simple $\Gamma$-representation such that $\alpha \geq \beta = (3,3;2,2,2)$, then $\tau \not= id$ on $\wis{iss}_{\alpha} \Gamma$ and there are simple representations $\phi \in \wis{iss}_{\alpha} \Gamma$ such that $Tr_{\phi}(b) \not= Tr_{\phi}(b')$.
\end{proposition}

\begin{proof} The unique open stratum of $\wis{iss}_{\alpha} \Gamma$ corresponds to the unique maximal representation type $\tau_{gen} = (1,\alpha)$, that is, $\wis{iss}_{\alpha} \Gamma(\tau_{gen})$ is the open set of simple $\Gamma$-representations. 

If $\alpha-\beta$ is the dimension vector of a simple $\Gamma$-representation, then we have a representation type $\tau = (1,\beta;1,\alpha-\beta)$ such that $\tau \not= id$ and $Tr(b) \not= Tr(b')$ on $\wis{iss}_{\alpha} \Gamma(\tau)$. But then, by the previous lemma, these facts also hold for $\wis{iss}_{\alpha} \Gamma(\tau_{gen})$.

If $\alpha-\beta$ is not the dimension vector of a simple $\Gamma$-representation, we consider the generic (maximal) representation type $\tau'=(e_1,\beta_1;\hdots;e_k,\beta_k)$ in $\wis{iss}_{\alpha-\beta} \Gamma$. But then, $\tau = (1,\beta;e_1,\beta_1;\hdots;e_k,\beta_k)$ is a representation type for $\wis{iss}_{\alpha} \Gamma$ and we can repeat the foregoing argument.
\end{proof}

\begin{proposition} If $\alpha=(a,b;x,y,z)$ is a simple dimension vector such that $\tau$ acts trivially on $\wis{iss}_{\alpha} B_3$, then
\[
dim~\wis{iss}_{\alpha} B_3 = n = a+b = x+y+z \]
\end{proposition}

\begin{proof} By the previous result we must have $\beta \not\leq \alpha$ and hence either $n \leq 5$ or $min(x,y,z)=1$. For a simple $B_3$-dimension vector we may assume that $a \geq b$ and $x=max(x,y,z)$, which leaves us with the following list of potential fixed-point components 
\[
\begin{array}{|c|l|c|}
\hline
n & \alpha & dim~\wis{iss}_{\alpha} B_3 \\
\hline
1 & (1,0;1,0,0) & 1 \\
2 & (1,1;1,1,0) & 2 \\
& (1,1;1,0,1) & 2 \\
3 & (2,1;1,1,1) & 3 \\
4 & (2,2;2,1,1) & 4 \\
5 & (3,2;2,2,1) & 5 \\
& (3,2;2,1,2) & 5 \\
\hline
6 & (3,3;3,2,1) & 6 \\
& (3,3;3,1,2) & 6 \\
& (4,2;2,2,2) & 6 \\
\hline
2k & (k,k;k,k-1,1) & 2k \\
& (k,k;k,1,k-1) & 2k \\
\hline
2k+1 & (k+1,k;k,k,1) & 2k+1 \\
& (k+1,k;k,1,k) & 2k+1 \\
\hline
\end{array}
\]
\end{proof}

By example~\ref{example} we know  that the $1$- and $2$-dimensional components are fixed-point components. All other potential fixed-point components belong to the infinite families, with one  exception: $(4,2;2,2,2)$. In the following sections we will prove that all of these are indeed fixed-point components.

\section{The infinite families}

In this section we will prove that  for $\alpha=(k,k;k,k-1,1)$ (the even case) and $\alpha=(k+1,k;k,k,1)$ (the odd case), $\wis{iss}_{\alpha} B_3$ is a fixed-point component. We will prove the even case by direct matrix calculations and deduce the odd case from it by a degeneration argument.

\begin{proposition} \label{infinite1} For all $k \in \mathbb{N}_+$ and $\alpha = (k,k;k,k-1,1)$, $\wis{iss}_{\alpha} B_3$ is a fixed-point component.
\end{proposition}

\begin{proof}
A general representation in $\wis{iss}_{\alpha} \Gamma$ corresponds to an  invertible $2m \times 2m$ matrix $B$ and quiver-representation
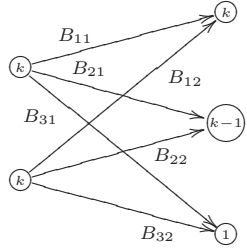
\[
\xymatrix@=.4cm{
& & & & \vtx{k} \\
\vtx{k} \ar[rrrru]^(.3){B_{11}} \ar[rrrrd]^(.3){B_{21}} \ar[rrrrddd]_(.2){B_{31}} & & & & \\
& & & & \vtx{{\tiny k-1}} \\
\vtx{k} \ar[rrrruuu]_(.7){B_{12}} \ar[rrrru]_(.7){B_{22}} \ar[rrrrd]_(.7){B_{32}} & & & & \\
& & & & \vtx{1}}
\]
After a base-change in the lower-left hand vertex, we may assume that the modified matrix-blocks become
\[
\begin{bmatrix} B'_{22} \\ B'_{32} \end{bmatrix} = I_k \]
The block $B_{12}$ is modified to an invertible $k \times k$ matrix $B'_{12}$ which becomes the identity matrix $I_k$ after performing a base-change in the top-right hand vertex. This changes the block $B_{11}$ to an invertible matrix $B'_{11}$ which becomes the identity matrix $I_k$ after a base-change in the top-left hand vertex. Hence, we may assume that, up-to isomorphism, the matrix $B$ has the following block-form
\[
B = \begin{bmatrix} B_{11} & B_{12} \\ B_{21} & B_{22} \\ B_{31} & B_{32} \end{bmatrix} = \begin{bmatrix} I_k & I_k \\ A & I_k \end{bmatrix} \]
with all blocks of sizes $k \times k$ and $A$ an invertible matrix such that $B$ is invertible.
One verifies that 
\[
(B^{-1})^{tr} = \begin{bmatrix} -C & I_k+C \\ C & -C \end{bmatrix}~\qquad \text{with}~\qquad C = (A-I_k)^{-1} \]
and performing the base-change
\[
(AC^{-1},-C;-A^{-1},I_{k-1},I_1) \in GL_k \times GL_k \times GL_k \times GL_{k-1} \times GL_1 \]
we obtain
\[
B = \begin{bmatrix} I_k & I_k \\ A & I_k \end{bmatrix} = \begin{bmatrix} -A^{-1} & 0 \\ 0 & I_k \end{bmatrix} \begin{bmatrix} -C & I_k + C \\ C & -C \end{bmatrix} \begin{bmatrix} C^{-1} A & 0 \\ 0 & -C^{-1} \end{bmatrix} \]
Therefore, the $\Gamma$-representations determined by the matrices $B$ and $(B^{-1})^{tr}$ are isomorphic and hence the involution $\tau$ is the identity map on the component $\wis{iss}_{\alpha} B_3$.
\end{proof}

\begin{proposition} \label{infinite2} For all $k \in \mathbb{N}_+$ and $\alpha = (k+1,k;k,k,1)$, $\wis{iss}_{\alpha} B_3$ is a fixed-point component.
\end{proposition}

\begin{proof}
Let $\alpha_+=(k+1,k+1;k+1,k,1)$, then the stratum
$\tau = (1,\alpha;1,(0,1;1,0,0))$ lies in the closure of the generic stratum $\tau_{gen}=(1,\alpha_+)$ in $\wis{iss}_{\alpha_+} \Gamma$. The result follows from the foregoing proposition and lemma~\ref{strata}.\end{proof}

\section{The exceptional component and vectorbundles on $\mathbb{P}_2$}

To finish the proof of theorem~\ref{main}, it suffices to show that $\wis{iss}_{\beta} B_3$ is a fixed-point component for $\beta=(4,2;2,2,2)$.  In \cite{LeBruynB3} we have given a parametrization of the matrices $B$ for a dense open subset of $\wis{iss}_{\beta} \Gamma$
\[
B=\begin{bmatrix} 1 & 0 & 0 & 0 & a & 0 \\
0 & 1 & e & 1 & 0 & 1 \\
1 & c & d & 0 & 1 & 0 \\
0 & 0 & 0 & 1 & 0 & b \\
0 & 1 & 0 & 0 & 1 & 0 \\
0 & 0 & 1 & 0 & 0 & 1 \end{bmatrix}
\]
One can attempt to show that $B$ and $(B^{-1})^{tr}$ belong to the same $GL(\beta)$-orbit by explicit computation. We follow a different approach, allowing us to connect this component to the study of stable vectorbundles on $\mathbb{P}_2$.

\begin{proposition} For $\alpha = (2n,n;n,n,n)$, the component $\wis{iss}_{\alpha} \Gamma$ is birational to $M_{\mathbb{P}^2}(n;0,n)$, the moduli space of semi-stable rank $n$ bundles on $\mathbb{P}_2$ with Chern classes $c_1=0$ and $c_2=n$.
\end{proposition}

\begin{proof} A representation in $\wis{rep}_{\alpha} \Gamma$ in general position
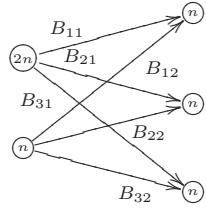
\[
\xymatrix@=.3cm{
& & & & \vtx{n} \\
\vtx{2n} \ar[rrrru]^(.3){B_{11}} \ar[rrrrd]^(.3){B_{21}} \ar[rrrrddd]_(.2){B_{31}} & & & & \\
& & & & \vtx{n} \\
\vtx{n} \ar[rrrruuu]_(.7){B_{12}} \ar[rrrru]_(.7){B_{22}} \ar[rrrrd]_(.7){B_{32}} & & & & \\
& & & & \vtx{n}}
\]
is such that  $\psi : \C^{2n} \rTo^{B_{11} \oplus B_{21} \oplus B_{31}} \C^n \oplus \C^n \oplus \C^n$ in injective, whence its cokernel defines maps
$
Cok(\psi)~:~\C^n \oplus \C^n \oplus \C^n \rTo^{(C_{12},C_{22},C_{32})} \C^n$
and therefore a representation for the quiver-setting
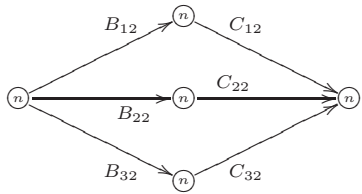
\[
\xymatrix@=.3cm{
& & & & \vtx{n} \ar[rrrrdd]^(.3){C_{12}} & & & & \\
 & & & & & & & & & \\
\vtx{n}  \ar[rrrruu]^(.7){B_{12}} \ar[rrrr]_(.7){B_{22}} \ar[rrrrdd]_(.7){B_{32}} & & & &   \vtx{n} \ar[rrrr]^(.3){C_{22}} & & & & \vtx{n} \\
 & & & & & & &  \\
& & & & \vtx{n} \ar[rrrruu]_(.3){C_{32}} & & & &  }
\]
By the general theory of reflection functors, isomorphism classes of representations are preserved under this construction. By the fundamental theorem of $GL_n$-invariants \cite[Thm. II.4.1]{Kraftbook} we can eliminate the base-change action in the middle vertices and obtain a representation of the quiver-setting
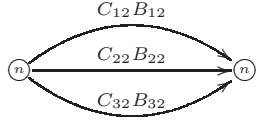
\[
\xymatrix@=.5cm{
\vtx{n} \ar@/^4ex/[rrrr]^{C_{12}B_{12}} \ar[rrrr]^{C_{22}B_{22}} \ar@/_4ex/[rrrr]^{C_{32}B_{32}} & & & & \vtx{n}}
\]
By results of  Klaus Hulek \cite{Hulek}, the corresponding moduli space of semi-stable quiver-representations (as in \cite{King} for the stability structure $(-1,1)$) is birational to $M_{\mathbb{P}_2}(n;0,n)$.
\end{proof}

\begin{proposition} $\wis{iss}_{\beta} B_3$ is a fixed point component.
\end{proposition}

\begin{proof} By results of Wolf Barth \cite{Barth}, we know that a stable rank $2$ bundle $\mathcal{E}$ on the projective plane with Chern-classes $c_1=0$ and $c_2=2$ is determined up to isomorphism by its curve of jumping lines, that is the collection of those lines $L \subset \mathbb{P}_2$ such that $\mathcal{E} | L \not\simeq \Oscr_L^{\oplus 2}$. If $\mathcal{E}$ is determined by the quiver-representation as in the previous proposition and if $x,y,z$ are projective coordinates of the dual plane $\mathbb{P}_2^*$, then the equation of this curve of jumping lines is given by
\[
det(C_{12}B_{12}x + C_{22}B_{22} y + C_{32}B_{32} z) = 0 \]
In terms of the matrix $B$ and its inverse $B^{-1}$ these $2 \times 2$ matrices are given as
\[
\underbrace{\begin{bmatrix} \ast & \ast & \ast \\ \ast & \ast & \ast \\ C_{12} & C_{22} & C_{32} \end{bmatrix}}_{B^{-1}}  \underbrace{\begin{bmatrix} \ast & \ast & B_{12} \\ \ast & \ast & B_{22} \\ \ast & \ast & B_{32} \end{bmatrix}}_{B} = \begin{bmatrix} I_2 & 0 & 0 \\ 0 & I_2 & 0 \\ 0 & 0 & I_2 \end{bmatrix} \]
But then, the bundle $\mathcal{F}$ corresponding to the matrix $(B^{-1})^{tr}$ is determined by the $2 \times 2$ matrices
\[
\underbrace{\begin{bmatrix} \ast & \ast & \ast \\ \ast & \ast & \ast \\ B_{12}^{tr} & B_{22}^{tr} & B_{32}^{tr} \end{bmatrix}}_{B^{tr}}  \underbrace{\begin{bmatrix} \ast & \ast & C_{12}^{tr} \\ \ast & \ast & C_{22}^{tr} \\ \ast & \ast & C_{32}^{tr} \end{bmatrix}}_{(B^{-1})^{tr}} = \begin{bmatrix} I_2 & 0 & 0 \\ 0 & I_2 & 0 \\ 0 & 0 & I_2 \end{bmatrix}  \]
and hence its curve of jumping lines
\[
det(B_{12}^{tr}C_{12}^{tr} x + B_{22}^{tr} C_{22}^{tr} y + B_{32}^{tr} C_{32}^{tr} z) \]
is the same as that for $\mathcal{E}$ and hence by Barth's result $\mathcal{E} \simeq \mathcal{F}$.
\end{proof}

\begin{remark} One can repeat the above argument verbatim for $\alpha = (2n,n;n,n,n)$.  However, if $n > 2$, the bundle $\mathcal{E}$ corresponding to the matrix $B$ is determined by its curve of jumping lines (defined as above by the $n \times n$ matrices $B_{ij}$ and $C_{ij}$) together with a half-canonical divisor on it, see \cite{Hulek}. Whereas the curve of jumping lines $Y$ of the bundle $\mathcal{F}$ corresponding to the matrix $(B^{-1})^{tr}$  coincides with that of $\mathcal{E}$, the involution $\tau$ acts non-trivially on the Jacobian $Pic^d_Y$ where $d=\frac{1}{2}n(n-1)$.
\end{remark}

\end{document}